\documentclass[12pt]{article}%
\usepackage{amsmath}
\usepackage{amsfonts}
\usepackage{amssymb}
\usepackage{graphicx}%
\providecommand{\U}[1]{\protect\rule{.1in}{.1in}}
\newtheorem{theorem}{Theorem}

\newtheorem{condition}{Condition}

\newtheorem{lemma}{Lemma}

\newtheorem{proposition}{Proposition}
\newtheorem{remark}{Remark}

\begin{document}

\title{High excursion probabilities for Gaussian fields on smooth manifolds. }
\author{Vladimir I. Piterbarg\thanks{Lomonosov Moscow state university, Moscow,
Russia; Scientific Research Institute of System Development of the Russian
Academy of sciences; International Laboratory of Stochastic Analysis and its
Applications, National Research University Higher School of Economics, Russian
Federation.\ \emph{{piter@mech.math.msu.su}}}}
\maketitle

\

\textbf{Abstract:}\footnotetext[1]{Partially funded by RFBR and CNRS, project
number 20-51-15001} Gaussian random fields on finite dimensional smooth
manifolds whose variances reach their maximum values at smooth submanifolds
are considered. Exact asymptotic behaviors of large excursion probabilities
have been evaluated. Vector Gaussian processes, chi-square processes, Bessel
process, fractional Bessel process, Bessel bridge are examples of application
of this result.

\

\textbf{Keywords:} Non-stationary random field; Gaussian vector process;
Gaussian field; large excursion; Pickands' method; Double sum method.

\section{ Introduction.}

This work is a continuation of papers \cite{KP}, \cite{KPR} and \cite{KHP},
where asymptotic behavior of large excursuion probabilities were evaluated for
Gaussian nonhomogeneous fields given on Euclidean space. An essential
condition in these and many previous works is that the variance of a
considered Gaussian field reaches its maximum at unique point in the
considered parameter set. In the mentioned three papers, general conditions
where given on behavior the variance near its maximum point. It has been shown
there that whereas introduced first in 1978, \cite{PP}, regular behavior of
the correlation function near the maximum variance point of the variance is
crucial for all known general tools for the asymptotic behavior evaluation,
introduced there behavior of the variance near this point can be much more general.

From the other hand, in many situations one have to consider Gaussian fields
with a set of maximum variance points. A typical situation is the large
excursions probability evaluation for norms of Gaussian vector processes, when
one uses the duality and passes to Gaussian fields on a cylinder, and theirs
variance can reach theirs maximum at a submanifold of the cylinder. Well known
examples are $\chi^{2}$ process and its generalizations, see \cite{book} for
examples and references, and Bessel process and bridge, with generalizations
to corresponding fracttion processes, \cite{PR}. These examples motivated the
present work.

\subsection{General setting}

Let $X(\mathbf{t}),$ $\mathbf{t}\in\mathcal{S\subset}\mathbb{R}^{n},$ be a
zero mean a.s. continuous Gaussian random field given on a $d$-dimensional
closed differentiable manifold $\mathcal{S},$ $d<n.$ Denote the covariance
function of $X$ by $R(\mathbf{s,t)=}EX(\mathbf{s})X(\mathbf{t})$, and by
$\sigma^{2}(\mathbf{t})=R(\mathbf{t,t}),$ the variance function of $X.$ Assume
the following,

\begin{condition}
\label{CondA2} Maximum points set of $\sigma(\mathbf{t}),$ $\mathbf{t\in}$
$\mathcal{S},$
\[
\mathcal{M}:\mathcal{=}\operatorname{argmax}_{\mathbf{t}\in\mathcal{S}}%
\sigma(\mathbf{t)}%
\]
is a differentiable finitely connected submanifold of $\mathcal{S}$ of
dimension $r<d.$ In particular, $r=0$ means that $\mathcal{M}$ consists of
finite number points on $\mathcal{S}.$
\end{condition}

\begin{remark}
\label{r=0}The case $r=0$ is considered in \cite{KP}, see also \cite{KPR}, so
that we assume further that $r\geq1.$The case $r=d,$ in particular
$\sigma(\mathbf{t)}$ is a constant on $\mathcal{S},$ is called a locally
homogeneous (stationary) case, provided the below Condition \ref{CondA3} is
fulfilled. We shall consider this case as well.
\end{remark}

Without loss of generality assume that $\sigma(\mathbf{t})\equiv1,$
$\mathbf{t\in}\mathcal{M}.$

We study the asymptotic behavior of the probability
\begin{equation}
P(\mathcal{S};u):=P(\max_{\mathbf{t\in}\mathcal{S}}X(\mathbf{t)>}u)
\label{P(S,u)}%
\end{equation}
as $u\rightarrow\infty.$ Remark immediately, that we shall use analogous
notation, $P(A;u),$ for any $A,$ meaning closure of $A$ if it is not close.
Remind that $\mathcal{S}$ is bounded.

\subsection{Examples. $\chi^{2}$ processes. Bessel process.}

Let $\mathbf{X}(t)=(X_{1}(t),...,X_{d}(t)),$ $t\in\lbrack0,T],$ be a Gaussian
vector process, $||\cdot||$ be a norm in some $d$-dimensional linear space
$\mathbb{L}^{d}$, generated by a scalar product $\left\langle \cdot
,\cdot\right\rangle .$ Denote by $\mathbb{M}^{d}$ the corresponding dual
space. By duality,
\[
P(\max_{t\in\lbrack0,T]}||\mathbf{X(}t)||>u)=P(\max_{(t,\mathbf{v)\in}T\times
S_{d}}\left\langle \mathbf{v,X}(t)\right\rangle >u),
\]
with $S_{d},$ the unit sphere in $\mathbb{M}^{d}.$

\textbf{Example. }$\chi^{2}$\textbf{ processes.} Let $X_{i}(t),$ $t\in
\lbrack0,T],$ $i=1,...,d,$ be independent Gaussian zero mean processes
$EX_{i}^{2}(t)=\sigma_{i}^{2}(t)$, consider generalized $\chi^{2}$ process,
\[
\chi_{\mathbf{b}}^{2}(t)=\sum_{j=1}^{d}b_{j}^{2}X_{j}^{2}(t),\ t\in
\lbrack0,T].
\]
If all $b_{j}$ are equal, we have usual $\chi^{2}$ process, the corresponding
field $\left\langle \mathbf{v,X}(t)\right\rangle $ is homogeneous with respect
to turns and shifts, if not, the variance of the field is equal to
\[
\sigma^{2}(t,\mathbf{b,v})=\sum_{j=1}^{d}\sigma_{i}^{2}(t)b_{j}^{2}v_{j}%
^{2},\mathbf{b}=(b_{1},...,b_{d}),
\]
with various structures of the set of maximum points.

For example, in case all $\sigma_{i}^{2}(t)$ are equal to $\sigma^{2}(t),$
with unique maximum point, and $b_{1}>b_{j},$ $j>1,$ the set $\mathcal{M}$
consists of two points, $\max\sigma^{2}(t,\mathbf{b,v})=b_{1}^{2}$; in case
$\sigma^{2}(t)\equiv1$ with the same $b_{j}$s, we have $\mathcal{M=}S_{d},$
\cite{chi}.

\textbf{Example. Bessel process.} For $X_{i}(t)=W_{i}(t),$ $t\in\lbrack0,1],$
i.i.d. Brownian motions, $\beta(t)=||\mathbf{X(}t)||$ is Bessel process. We
have,%
\[
P(\max_{[0,1]}\beta(t)>u)=P(\max_{(t,\mathbf{v)\in}T\times S_{d}}\sum
_{i=1}^{d}v_{i}W_{i}(t)>u).
\]
Here $\mathcal{M}=\{1\}\times S_{d}.$ Using behavior of the covariance
function of $\mathbf{X(}t)$ near this set, it is shown in \cite{PR} that
\[
P(\max_{[0,1]}\beta(t)>u)=\frac{\pi^{(d-1)/2}}{2^{d/2-1}\Gamma(d/2)}%
u^{d-2}e^{-u^{2}/2}(1+o(1)),
\]
as $u\rightarrow\infty.$ Also Bessel bridge with zero end point, fractional
Bessel processes and bridges have been similarly considered in \cite{PR},
definitions are clear.

\subsection{Tools. Pickands' Double Sum Method.}

An idea how to consider non-smooth stationary Gaussian processes belong to J.
Pickands III, 1969. He first considered the probability on infinitely small
intervals and then passed to arbitrary interval, using semi-additivity of the
probability with Bonferroni inequalities. It is why it is called Double Sun
Method. A generalization to non-stationary processes has been done in
\cite{PP}. Let $X(t),$ $t\in\lbrack0,1]$ be a zero mean Gaussian process with
correlation function $r(s,t)$ and variance $\sigma^{2}(t),$ which has a unique
maximum point $t_{0}\in(0,1).$ Assume for some positive $a,\beta$ and
$\alpha\leq2$ that
\[
\sigma(t)=1-a|t-t_{0}|^{\beta}(1+o(1)),\ \text{\ \ }t\rightarrow t_{0};
\]%
\[
r(s,t)=1-|t-s|^{\alpha}(1+o(1)),\ \text{\ \ }s\rightarrow t_{0},\ t\rightarrow
t_{0};
\]
and for some positive $\gamma,G$ and all $s,t,$
\[
E(X(t)-X(s))^{2}\leq G|t-s|^{\gamma}.
\]

Denote
\[
\Psi(u)=\frac{1}{\sqrt{2\pi}u}e^{-u^{2}/2},
\]
the asymptotic of the tail of Gaussian standard distribution. In the above
conditions, for some constants $C_{1}>0$ and $C_{2}>1,$ see Lemma
\ref{lemma_pickands} below, we have as $u\rightarrow\infty$ the following.

If $\beta>\alpha$ (\emph{stationary like case}),
\[
P([0,1],u)=C_{1}u^{2/\alpha-2/\beta}\Psi(u)(1+o(1)).
\]

If $\beta=\alpha$ (\emph{transition case}),
\[
P([0,1],u)=C_{2}\Psi(u)(1+o(1)).
\]

If $\beta<\alpha$ (\emph{Talagrand case}),
\[
P([0,1],u)=\Psi(u)(1+o(1)).
\]

From the proof of assertions below it will be cleare some physical sence of
these asymptotic relations.

In the stationary like case, since the correlation is sharper at $t_{0}$ than
variance, trajectories oscillate near $t_{0},$ do not paying much attention on
slow behavior of the variance, therefore we have a degree of $u$ before
$\Psi.$

In the Talagrand case, conversely, the variance does not pay any attention on
the slow trajectories oscillations, M. Talagrand, \cite{talagrand}, showed
this relation in the maximal generality.

In transition case, the behaviors of variance and correlation are equal one to
another, so the order of asymptotic behavior is the same, but with constant
greater than one.

It turned out, that, while power behavior of the correlation is critical for
Pickands' method, behavior of the variance can be much more general. We shall
see, that the only important point is whether the variance behavior is slower
or faster, or equivalent the behavior of correlation. In case of processes, it
was described in \cite{KHP}. Then in \cite{KP}, with corrections in
\cite{KPR}, this have been generalized on Gaussian fields with single maximum
variance point$.$

The present work is an immediate continuation of these works, several
assertions here have almost the same proofs, therefore they are mainly
omitted. The only constructions are given in details which are significantly
different in case $\dim\mathcal{M}>0$ than the corresponding ones in case
$\dim\mathcal{M}=0.$

\section{Extracting of an informative set.}

Our first aim is to extract an informative small neighbourhood of
$\mathcal{M}$ in $\mathcal{S}$ such that the probability (\ref{P(S,u)}) is
equivalent to the same probability but with this neighbourhood instead of
$\mathcal{S}.$

Denote by $\mathcal{M}_{\varepsilon}\subset\mathcal{S},$ an $\varepsilon
$-neighborhood of $\mathcal{M}$ in $\mathcal{S\subset}$ $\mathbb{R}^{n},$ that
is,%
\begin{equation}
\mathcal{M}_{\varepsilon}:=\left\{  \mathbf{s}:\min_{\mathbf{t}\in\mathcal{M}%
}\mathbf{||t}-\mathbf{s||<}{\varepsilon}\right\}  . \label{M_eps}%
\end{equation}
It is easy to show, using Borell-Tsirelson-Ibragimov-Sudakov (Borell-TIS)
inequality, that for any $\varepsilon>0$ there exists $\delta>0$ such that
\begin{equation}
P(\mathcal{S};u)=P(\max_{\mathbf{t\in}\mathcal{M}_{\varepsilon}}%
X(\mathbf{t)>}u)(1+O(e^{-\delta u^{2}}),\ \ u\rightarrow\infty, \label{extr0}%
\end{equation}
see, for example, Theorem D.1, \cite{book}.

Further, we need in proofs a more narrow neighbourhood of $\mathcal{M}.$ To
this end, assume a slightly stronger condition than a.s. continuity of sample paths.

\begin{condition}
\label{CondA1} There exists $\varepsilon>0$ such that Dudley's entropy
integral with respect to pseudo-semi-metrics generated by the standardized
field $\bar{X}(\mathbf{t})=X(\mathbf{t})/\sigma(\mathbf{t})$, $\mathbf{t}%
\in\mathcal{M}_{\varepsilon},$ is finite. For definitions see \cite{Dudley},
\cite{lectures}, \cite{book}.
\end{condition}

Notice that for homogeneous Gaussian fields this condition is also necessary
for existing of a.s. continuous version of the field (X. Fernique
\cite{fernique}).

From this condition it follows V. A. Dmitrovsky's inequality, \cite{Dmitr1},
\cite{Dmitr2}, \cite{lectures}, which is more accurate than Borell-TIS one.
Using this inequality, it is proved in \cite{KP} the following Lemma.

\begin{lemma}
\label{extracting} In the above conditions and notations, there exists
$\gamma(u)$ with $\gamma(u)\rightarrow0$ as $u\rightarrow\infty$, such that
\begin{equation}
P(\mathcal{S};u)=P(\mathcal{M}_{u};u)\left(  1+O\left(  e^{-\log^{2}u}\right)
\right)  \label{ext1}%
\end{equation}
as $u\rightarrow\infty,$ where
\begin{equation}
\ \ \mathcal{M}_{u}:=\mathcal{M}_{\varepsilon(u)}\ \text{with }\varepsilon
(u)={2}u^{-1}\gamma(u)+2u^{-2}\log^{2}u. \label{M_e}%
\end{equation}

\end{lemma}

Thus we are in a position to study the asymptotic behavior of the latter probability.

\section{Exact asymptotic behavior.}

Denote%
\[
r(\mathbf{s}_{1}\mathbf{,s}_{2}\mathbf{)=}\frac{R(\mathbf{s}_{1}%
\mathbf{,s}_{2}\mathbf{)}}{\sigma(\mathbf{s}_{1})\sigma(\mathbf{s}%
_{2}\mathbf{)}},
\]
the correlation function of $X.$

\begin{condition}
\label{CondA3} (Local homogeneity). There exists $\varepsilon>0$ such that for
any $\mathbf{t}\in\mathcal{M}_{\varepsilon}$ there exists a correlation
function $r_{\mathbf{t}}(\mathbf{s),}$ $\mathbf{s}\in\mathbb{R}^{n}$ of a
homogeneous random field with $r_{\mathbf{t}}(\mathbf{s)}<1$ for all
$\mathbf{s\neq0}$ and such that
\[
\lim_{\mathbf{s}_{1}\mathbf{,s}_{2}\rightarrow\mathbf{t}}\frac{1-r(\mathbf{s}%
_{1}\mathbf{,s}_{2}\mathbf{)}}{1-r_{\mathbf{t}}(\mathbf{s}_{2}\mathbf{-s}%
_{1})}=1.
\]

\end{condition}

\textbf{The atlas on }$\mathcal{M}_{\varepsilon}.$ Remind some definitions
related to the smooth manifold $\mathcal{M}_{\varepsilon}$, see for example
\cite{zorich}. Any $\mathbf{t\in}\mathcal{M}$ has a neighbourhood
$\mathcal{U}\subset$ $\mathcal{S}$ which is homeomorphic to $\mathbb{R}^{d},$
and the corresponding homeomorphism $\varphi:\mathbb{R}^{d}\rightarrow
\mathcal{U}$ is differentiable. We call the set $U=\varphi^{\leftarrow
}(\mathcal{U}),$ \textit{the map of }$\mathcal{U}.$ Since $\mathcal{M}%
_{\varepsilon}$ is closed, there exists a finite collection of such
neighbourhoods which cover whole $\mathcal{M}$ and such that for all
sufficiently small $\varepsilon$ the union of them cover whole $\mathcal{M}%
_{\varepsilon}.$ Such the collection together with corresponding maps forms an
\emph{atlas} on $\mathcal{M}_{\varepsilon},$ which we fix for future considerations.

\textbf{Agreement on coordinates. }Now consider $\mathcal{U}$ from the fixed
atlas, with corresponding $\varphi,$ and choose the coordinate system in
$\mathbb{R}^{n}$ such that in the coordinates of $\mathbf{s}=(s_{i}%
,i=1,...,n)\mathbf{\in}\mathbb{R}^{n},$
\begin{align}
U  &  :=\varphi^{\leftarrow}(\mathcal{U)\subset}\mathbb{R}^{d}:=\{\mathbf{s\in
}\mathbb{R}^{n}:\{s_{d+1}=...=s_{n}=0\}\ \nonumber\\
&  \ \text{and }\nonumber\\
M  &  :=\varphi^{\leftarrow}(\mathcal{U\cap M)\subset}\mathbb{R}%
^{r}:=\{\mathbf{s\in}\mathbb{R}^{n}:\{s_{r+1}=...=s_{n}=0\}.
\label{local coord}%
\end{align}
From now on we use these coordinates and write $\mathbf{s=}(s_{1}%
,...,s_{d})\in U.$ For such $\mathbf{s}$ denote
\begin{equation}
\mathbf{s}^{1}=(s_{1},...,s_{r},0,...,0)\in M, \label{ort1}%
\end{equation}
the projection of $\mathbf{s}$ on $M,$ and
\begin{equation}
\mathbf{s}^{2}=(0,...,0,s_{r+1},...,s_{d}), \label{ort2}%
\end{equation}
so that $\mathbf{s=s}^{1}+\mathbf{s}^{2},$ the orthogonal decomposition of
$\mathbf{s.}$ We will write sometimes $\mathbf{s=(s}^{1},\mathbf{s}^{2}).$

Thus for any map from the fixed atlas,%
\begin{equation}
P(\mathcal{U};u)=P(\max_{\mathbf{s\in}U}X(\mathbf{s)>}u)=P(U;u), \label{PUu}%
\end{equation}
where $\mathbf{s}\in\mathbb{R}^{d}$ is scripted in the selected coordinates.
As well, with (\ref{extr0}, \ref{ext1}) we have in mind the same relation s
with $M_{\varepsilon}$ and $M_{u}$ instead of $\mathcal{M}_{\varepsilon}$ and
$\mathcal{M}_{u},$ correspondingly.

Notice immediately, that the notation $d\mathbf{t}$ is used in integrals below
for volume elements in $U,$ $\mathcal{U}$ and theirs subsets, do not depending
of dimensions of the integral sets. Transitions from one integral to the
corresponding another one perform in standard way.

We find first the asymptotic behavior of probability of $P(M_{\varepsilon
};u),$ and then, using introduced Conditions and the fact that the atlas is
finite, pass to the probability (\ref{P(S,u)}). Write for short $M$ instead of
$M\cap U.$

For points $\mathbf{a}=(a_{1},...,a_{d}),$ $\mathbf{b}=(b_{1},...,b_{d})$ in
$S$ denote $\mathbf{ab}=(a_{1}b_{1},...,a_{d}b_{d})\in\mathbb{R}^{d}.$ For a
set $T\subset\mathbb{R}^{d},$ we write $\mathbf{a}T=\{\mathbf{at,t\in}T\}.$
The cases $\mathbf{ab\in}\mathcal{S},$ $T\subset\mathcal{S}$, $\mathbf{a}%
T\subset\mathcal{S}$ will be specified.

\begin{condition}
\label{CondA4} There exist $\varepsilon>0$, a vector function $\mathbf{q}%
(u)=(q_{1}(u),...,q_{d}(u)),$ $u>0,$ $q_{i}(u)>0,\ i=1,...,d$ and a positive
and finite for all $\mathbf{s\in}\mathcal{M}_{\varepsilon}\setminus
\mathcal{M}$ function $h(\mathbf{s})$ such that for any $\mathcal{U}$ from the
atlas, any $\mathbf{t}\in\mathcal{U\cap M},$ there exists on the corresponding
by (\ref{local coord}) maps $U,M$ a continuous vector function $\mathbf{C}%
_{\mathbf{t}}=(C_{1\mathbf{t}},...,C_{d\mathbf{t}}),$ $\mathbf{t}\in
M_{\varepsilon},$ and positive $c,C$ with $C_{i\mathbf{t}}\in\lbrack
c,C],\ i=1,...,d,$ such that for any $\mathbf{s}\in M_{\varepsilon}$ and for
$r_{\mathbf{t}}(\cdot)$ defined in Condition \ref{CondA3},
\begin{equation}
\lim_{u\rightarrow\infty}u^{2}(1-r_{\mathbf{t}}(\mathbf{C}_{\mathbf{t}%
}\mathbf{q}(u)\mathbf{s))=}h(\mathbf{s}) \label{A3}%
\end{equation}
uniformly in $\mathbf{s}$.
\end{condition}

Sometimes we shall write for short $\mathbf{q}_{\mathbf{t}}(u)=\mathbf{C}%
_{\mathbf{t}}\mathbf{q}(u)$ with corresponding components $q_{i\mathbf{t}%
}(u)=C_{i\mathbf{t}}q_{i}(u),$ $i=1,...,d.$

Remark that by definition of uniform convergence to a positive for all
non-zero $\mathbf{s}$ function with $h_{\mathbf{t}}(\mathbf{0})=0$, from
Condition \ref{CondA4} it follows that $h_{\mathbf{t}}(\mathbf{s})$ is
continuous, and for any continuous function $c_{\mathbf{s}}$ with
$\lim_{\mathbf{s\rightarrow0}}c_{\mathbf{s}}=1,$%
\begin{equation}
\lim_{\mathbf{s\rightarrow0}}\frac{1-r_{\mathbf{t}}(c_{\mathbf{s}}\mathbf{s)}%
}{1-r_{\mathbf{t}}(\mathbf{s})}=1\mathbf{.} \label{A4uniform}%
\end{equation}

Remind that we assumed that the basis in $\mathbb{R}^{d}$ satisfies Condition
\ref{CondA4}. From Conditions \ref{CondA3} and \ref{CondA4} it follows a
regular variation property of $r_{\mathbf{t}}(\mathbf{s})$. The following
proposition is proved in \cite{KP}.

\begin{proposition}
\label{e_reg_var}Let Conditions \ref{CondA3} and \ref{CondA4} be fulfilled for
a covariance function $r(\mathbf{t}).$ Then for any vector $\mathbf{f,}$
function $1-r_{\mathbf{t}}(s\mathbf{f)}$ regularly varies in $s$ at zero with
degree $\alpha(\mathbf{f})\in(0,2],$ and
\[
h(s\mathbf{f})=A_{\mathbf{f}}|s|^{\alpha(\mathbf{f})},\ \ A_{\mathbf{f}}>0.
\]
The function $\mathbf{q}(u)\mathbf{f}$ is regularly varying in $u$ at infinity
with degree $-2/\alpha(\mathbf{f}).$
\end{proposition}

It is easy to see that $\alpha(\mathbf{f})$ is equal to one of the
$\alpha(\mathbf{e}_{i}),$ $\{\mathbf{e}_{i},i=1,...,d\}$ is the basis of
$\mathbb{R}^{d}.$ As it is shown in \cite{KP}, if $\alpha(\mathbf{f})=2,$ then
$\lim_{t\rightarrow0}t^{-2}(1-r(t\mathbf{f}))>0,$ and for some $\ q_{-}>0,$%
\begin{equation}
\mathbf{q}(u)\mathbf{f}\geq q_{-}u^{-1},\ i=1,...,d. \label{q for Talagrand}%
\end{equation}

Now assume a behavior of $\sigma(\mathbf{t})$ near the points of $M.$ We shall
see from the proof of Lemma \ref{lemma_pickands} that the crucial point is the
behavior of the ratio
\[
\frac{1-\sigma(\mathbf{t}+\mathbf{C}_{\mathbf{t}}\mathbf{q}(u)\mathbf{s}%
)}{1-r_{\mathbf{t}}(\mathbf{C}_{\mathbf{t}}\mathbf{q}(u)\mathbf{s)}}%
\]
as $u\rightarrow\infty,$ where $\mathbf{t}\in M$ and $\mathbf{s}\in
M_{\varepsilon}$.

In view of Condition \ref{CondA4} we assume the following.

\begin{condition}
\label{CondA5} For any $\mathbf{t}\in\mathcal{M}$ and all $\mathbf{s}%
\in\mathcal{M}_{\varepsilon}$ there exists in the corresponding
$(M,M_{\varepsilon})$ the limit
\begin{equation}
h_{1\mathbf{t}}(\mathbf{s}):=\lim_{u\rightarrow\infty}u^{2}(1-\sigma
(\mathbf{t}+\mathbf{C}_{\mathbf{t}}\mathbf{q}(u)\mathbf{s})\in\lbrack
0,\infty]. \label{cond2}%
\end{equation}

\end{condition}

As in Subsection 1.3, in case when the limit is equal to zero we speak about
the \emph{ stationary like case}. If the limit is equal to infinity, we refer
to the \emph{Talagrand case}, since M. Talagrand, \cite{talagrand}, see
comments in \cite{KP} and \cite{book}. At last, we say about the \emph{
transition case }if $h_{1\mathbf{t}}(\mathbf{s})$ is neither zero nor infinity.

Denote correspondingly for any $\mathbf{t}\in M,$
\begin{align}
K_{0\mathbf{t}}  &  :=\{\mathbf{s}\in M_{\varepsilon}\setminus M\mathbf{:\ }%
h_{1\mathbf{t}}(\mathbf{s)=}0\},\ \ K_{c\mathbf{t}}:=\{\mathbf{s}\in
M_{\varepsilon}\setminus M\mathbf{:\ }h_{1\mathbf{t}}(\mathbf{s)\in(}%
0,\infty)\},\ \ \nonumber\\
K_{\infty\mathbf{t}}  &  :=\{\mathbf{s}\in M_{\varepsilon}\setminus
M\mathbf{:\ }h_{1\mathbf{t}}(\mathbf{s)=\infty}\}. \label{cases}%
\end{align}
We shall see that properties of union in $\mathbf{t}$ of these sets together
with all above conditions follow asymptotic behavior of the probability
$P(\mathcal{S};u)$.

\textbf{Pickands' Lemma. }Introduce a Gaussian a.s. continuous field
$\chi(\mathbf{s)}$ with $\chi(\mathbf{0)=}0,$ and
\[
\operatorname*{var}(\chi(\mathbf{s}_{1}\mathbf{)-}\chi(\mathbf{s}%
_{2}\mathbf{))}=2h(\mathbf{s}_{1}\mathbf{-s}_{2}\mathbf{),\ }E\chi
(\mathbf{s)}=-h(\mathbf{s})\mathbf{.}%
\]
The existence of such the field follows from Condition \ref{CondA4} and the
proof of Lemma below which, in fact, is a generalization of Lemma 6.1,
\cite{book}. For any $T\subset\mathbb{R}^{d}$ and any$\ \mathbf{t}\in M$
introduce the Pickands' type constant,
\begin{equation}
P_{\mathbf{q,t}}(T\mathbf{)=}E\exp(\max_{\mathbf{s\in}T\cap U}\chi
(\mathbf{s})-h_{1\mathbf{t}}(\mathbf{s)}).\ \mathbf{\ \ } \label{pit}%
\end{equation}

\begin{lemma}
\label{lemma_pickands}In the above notations and conditions, for any bounded
closed set $T\subset\mathbb{R}^{d}$ and any point $\mathbf{t}\in U\cap M,$
$\mathbf{t\notin}\partial U$%
\begin{equation}
P(\max_{\mathbf{s\in t+q}_{\mathbf{t}}(u)T}X(\mathbf{s)>}%
u)=(1+o(1))P_{\mathbf{q,t}}(T\mathbf{)}\Psi(u) \label{LP1}%
\end{equation}
as $u\rightarrow\infty,$ where $\mathbf{q}_{\mathbf{t}}(u)=\mathbf{C}%
_{\mathbf{t}}\mathbf{q}(u),$ as above. If $\mathbf{t}\in\partial U\cap M\ $and
$h_{1\mathbf{t}}(\mathbf{s)}\equiv0$ this relation is valid as well. If
$\mathbf{t}\in\partial U\cap M\ $and $h_{1\mathbf{t}}(\mathbf{s})\neq0$ for
some $\mathbf{s,}$ the set $T$ on the right will be correspondingly truncated,
denote it by $T_{\mathbf{t}}.$
\end{lemma}

\textbf{Outline of proof. }The proof mainly follows the proof of Lemma 6.1,
\cite{book}. Fix $\varepsilon$ in (\ref{M_eps}). Observe that since
$M_{\varepsilon}$ is open, $\mathbf{t+q}_{\mathbf{t}}(u)T\subset
M_{\varepsilon}$ for all sufficiently large $u.$ Further proof can be proceed
in $\mathbb{R}^{d}$ similarly to the Lemma 6.1 proof with obvious adition of
the case $\mathbf{t}\in M\cap\partial U.$ We repeat here initial evaluations
which are based on the conditions here. We have,%
\begin{align}
&  P(\max_{\mathbf{s\in t+q}_{\mathbf{t}}(u)T}X(\mathbf{s)>}u)\nonumber\\
&  =\frac{1}{\sqrt{2\pi}}\int_{-\infty}^{\infty}e^{-v^{2}/2}P\left(  \left.
\max_{\mathbf{s\in t+q}_{\mathbf{t}}(u)T}X(\mathbf{s)>}u\right\vert
X(\mathbf{t})=v\right)  dv\nonumber\\
&  =\frac{1}{\sqrt{2\pi}u}e^{-u^{2}/2}\int_{-\infty}^{\infty}e^{w-w^{2}%
/2u^{2}}P\left(  \left.  \max_{\mathbf{s\in t+q}_{\mathbf{t}}(u)T}%
X(\mathbf{s)>}u\right\vert X(\mathbf{t})=u-\frac{w}{u}\right)  dw,
\label{Pickands Lemma19}%
\end{align}
with change $v=u-w/u.$ Introduce the Gaussian process
\begin{equation}
\chi_{u\mathbf{t}}(\mathbf{s}):=u(X(\mathbf{t}+\mathbf{q}_{\mathbf{t}%
}(u)\mathbf{s})-u)+w. \label{chi_u}%
\end{equation}
The latter integral can be rewritten as
\begin{equation}
\int_{-\infty}^{\infty}e^{w-w^{2}/2u^{2}}P\left(  \left.  \max_{\mathbf{s}\in
T}\chi_{u\mathbf{t}}(\mathbf{s})>w\right\vert X(\mathbf{t})=u-\frac{w}%
{u}\right)  dw. \label{PL1}%
\end{equation}
Further, by (\ref{chi_u}), using formulas for conditional mean and variance,
\begin{align}
E  &  \left(  \left.  \chi_{u\mathbf{t}}(\mathbf{s})\right\vert X(\mathbf{t}%
)=u-\frac{w}{u}\right)  =\nonumber\\
&  =-u^{2}(1-R(\mathbf{t,t}+\mathbf{q}_{\mathbf{t}}(u)\mathbf{s}%
))+w(1-R(\mathbf{t,t}+\mathbf{q}_{\mathbf{t}}(u)\mathbf{s})), \label{dence09}%
\end{align}%
\begin{align}
&  \operatorname*{var}\left(  \left.  \chi_{u\mathbf{t}}(\mathbf{s}_{1}%
)-\chi_{u\mathbf{t}}(\mathbf{s}_{2})\right\vert X(\mathbf{t})=u-\frac{w}%
{u}\right) \nonumber\\
&  =u^{2}\left(  \operatorname*{var}[X(\mathbf{t}+\mathbf{q}_{\mathbf{t}%
}(u)\mathbf{s}_{1}))-X(\mathbf{t}+\mathbf{q}_{\mathbf{t}}(u)\mathbf{s}%
_{2}))]\right. \nonumber\\
&  \left.  \ \ \ \ \ \ \ \ \ \ \ \ \ \ \ -[R(\mathbf{t,t}+\mathbf{q}%
_{\mathbf{t}}(u)\mathbf{s}_{1})-R(\mathbf{t,t}+\mathbf{q}_{\mathbf{t}%
}(u)\mathbf{s}_{2})]^{2}\right)  . \label{dence9}%
\end{align}
Obviously,
\[
E\left(  \left.  \chi_{u\mathbf{t}}(\mathbf{0})\right\vert X(\mathbf{t}%
)=u-\frac{w}{u}\right)  =E\left(  \left.  \chi_{u\mathbf{t}}^{2}%
(\mathbf{0})\right\vert X(\mathbf{t})=u-\frac{w}{u}\right)  =0.
\]
Since $\sigma(\mathbf{t})=1,$ by Conditions \ref{CondA4} and \ref{CondA5},
\begin{align}
&  u^{2}(1-R(\mathbf{t,t}+\mathbf{q}_{\mathbf{t}}(u)\mathbf{s}))=u^{2}%
(1-r(\mathbf{t,t}+\mathbf{q}_{\mathbf{t}}(u)\mathbf{s})\sigma(\mathbf{t}%
+\mathbf{q}_{\mathbf{t}}(u)\mathbf{s))}\nonumber\\
&  =u^{2}(1-\sigma(\mathbf{t}+\mathbf{q}_{\mathbf{t}}(u)\mathbf{s}%
))+u^{2}(1-r(\mathbf{t,t}+\mathbf{C}_{\mathbf{t}}\mathbf{q}(u)\mathbf{s}%
))\sigma(\mathbf{t}+\mathbf{C}_{\mathbf{t}}\mathbf{q}(u)\mathbf{s}%
)\label{LP2proof}\\
&  =(h_{1\mathbf{t}}(\mathbf{s})+h(\mathbf{s}))(1+o(1)),\ \text{as
}u\rightarrow\infty.\nonumber
\end{align}
Further,
\begin{align*}
&  \operatorname*{var}[X(\mathbf{t}+\mathbf{q}_{\mathbf{t}}(u)\mathbf{s}%
_{1}))-X(\mathbf{t}+\mathbf{q}_{\mathbf{t}}(u)\mathbf{s}_{2}))]=\sigma
^{2}(\mathbf{t}+\mathbf{q}_{\mathbf{t}}(u)\mathbf{s}_{1})+\sigma
^{2}(\mathbf{t}+\mathbf{q}_{\mathbf{t}}(u)\mathbf{s}_{2})\\
&  -2\sigma(\mathbf{t}+\mathbf{q}_{\mathbf{t}}(u)\mathbf{s}_{1})\sigma
(\mathbf{t}+\mathbf{q}_{\mathbf{t}}(u)\mathbf{s}_{2})r(\mathbf{t}%
+\mathbf{q}_{\mathbf{t}}(u)\mathbf{s}_{1}),\mathbf{t}+\mathbf{q}_{\mathbf{t}%
}(u)\mathbf{s}_{2}))\\
&  =((\sigma(\mathbf{t}+\mathbf{q}_{\mathbf{t}}(u)\mathbf{s}_{1}%
)-1)+(1-\sigma(\mathbf{t}+\mathbf{q}_{\mathbf{t}}(u)\mathbf{s}_{2})))^{2}\\
&  +2\sigma(\mathbf{t}+\mathbf{q}_{\mathbf{t}}(u)\mathbf{s}_{1})\sigma
(\mathbf{t}+\mathbf{q}_{\mathbf{t}}(u)\mathbf{s}_{2})(1-r(\mathbf{t}%
+\mathbf{q}_{\mathbf{t}}(u)\mathbf{s}_{1}),\mathbf{t}+\mathbf{q}_{\mathbf{t}%
}(u)\mathbf{s}_{2})).
\end{align*}
By Condition \ref{CondA5}, the first summand on the right is equal to
$O(u^{-4})$ uniformly in $\mathbf{s}_{1},\mathbf{s}_{2}$ as $u\rightarrow
\infty.$ By conditions \ref{CondA3} and \ref{CondA4}, the second summand on
the right is equal to $2u^{-2}h(\mathbf{s}_{2}-\mathbf{s}_{1})(1+o(1))$
uniformly in $\mathbf{s}_{1},\mathbf{s}_{2}$ as $u\rightarrow\infty.$

The further proof is just a repetition of the Lemma 6.1 proof, including weak
and majorized convergence as $u\rightarrow\infty$ of the process
$\chi_{u\mathbf{t}}(\mathbf{\cdot})$ and the integral in
(\ref{Pickands Lemma19}).

\begin{remark}
\label{H_q}In the homogeneous like case, $h_{1\mathbf{t}}(\mathbf{s})\equiv0,$
we have,%
\[
P_{\mathbf{q,t}}(T\mathbf{)=:}H_{\mathbf{q}}(T\mathbf{)=}E\exp(\max
_{\mathbf{s\in}T}\chi(\mathbf{s})-h(\mathbf{s)}).\mathbf{\ \ }%
\]
Moreover, from the proof of standard Pickands theorem for homogeneous Gaussian
fields, see for example \cite{book} and a slight its extension in \cite{KP},
it follows that
\begin{equation}
\lim_{T\rightarrow\infty}T^{-d}H_{\mathbf{q}}([0,T]^{d})=:H_{\mathbf{q}}%
\in(0,\infty), \label{PP2_1}%
\end{equation}
the Pickands' constant. In Talagrand case, $h_{1\mathbf{t}}(\mathbf{s}%
)\equiv\infty,$ that is, $K_{\infty\mathbf{t}}\neq\emptyset$
\[
P_{\mathbf{q,t}}(T)=1.
\]

\end{remark}

\begin{remark}
\label{M_u} In fact, by Lemma \ref{extracting}, we need this lemma and the
following splitting and other construction only in the set $M_{u},$ but for
convenience we assume the informative set to be independed of $u,$ that is,
simply write $M_{\varepsilon}.$
\end{remark}

\textbf{Splitting. }In order to apply local Lemma \ref{lemma_pickands} we
construct a splitting of $M_{\varepsilon}$ into small sets. First split $M,$
and begin with the first axis of $\mathbb{R}^{d}$, $\{\mathbf{t=(}%
t_{1},0,...,0),t_{1}\in\mathbb{R}\}\subset M.$

Denote $\mathbf{T}^{1}:=(T,0,...,0)$ and take for $k_{1}=0,1,2,...,$
\[
\mathbf{t}_{0}^{1}=\mathbf{0},\mathbf{t}_{1}^{1}=\mathbf{q}_{\mathbf{t}%
_{0}^{1}}(u)\mathbf{T}^{1},\mathbf{t}_{2}^{1}=\mathbf{t}_{1}^{1}%
+\mathbf{q}_{\mathbf{t}_{1}^{1}}(u)\mathbf{T}^{1},...,\mathbf{t}_{k_{1}+1}%
^{1}=\mathbf{t}_{k_{1}}^{1}+\mathbf{q}_{\mathbf{t}_{k_{1}}^{1}}(u)\mathbf{T}%
^{1},
\]
till $\mathbf{t}_{k_{1}+1}^{1}$ such that $\mathbf{t}_{k_{1}}^{1}\in M,$ that
is the last point is out of $M.$ Similarly backward for $k_{1}=-1,-2,...,$till
$\mathbf{t}_{k_{1}-1}^{1}$ such that $\mathbf{t}_{k_{1}}^{1}\in M.$ Further,
at any point $\mathbf{t}_{k_{1}}^{1}=(t_{k_{1}1}^{1},0,...,0),$ construct a
similar grid on the axis $\{(t_{k_{1}1}^{1},t_{2},...,0),t_{2}\in
\mathbb{R}\}\subset M.$ Namely, denote $\mathbf{T}^{2}:=(0,T,0,...,0),$ and
put for $k_{2}=0,1,...$
\[
\mathbf{t}_{k_{1}0}^{2}=(t_{k_{1}1}^{1},0,...,0),...,\mathbf{t}_{k_{1}k_{2}%
+1}^{2}=(t_{k_{1}1}^{1},t_{k_{1}(k_{2}+1)1}^{2},0,...,0)=\mathbf{t}%
_{k_{1}k_{2}}^{2}+\mathbf{q}_{\mathbf{t}_{k_{1}k_{2}}^{2}}(u)\mathbf{T}%
^{2},...,
\]
till $\mathbf{t}_{k_{1}k_{2}+1}^{2}$ such that $\mathbf{t}_{k_{1}k_{2}}^{2}\in
M.$ Similarly backward for $k_{2}=-1,-2,...,$till $\mathbf{t}_{k_{1}k_{2}%
-1}^{2}\in M.$ So on, till $\mathbf{T}^{r}:=(0,...0,T,0,...,0),$ with $T$ on
$r$th place. Denote by $\mathcal{N}_{r}:=\{\mathbf{t}_{\mathbf{k}}^{r},$
$\mathbf{k\in}$ $\mathbb{Z}^{r}\},$ the constructed net on $M.$ Observe that
for any $\mathbf{t}_{\mathbf{k}}^{r},$ $2^{r}$ points $\mathbf{t}%
_{\mathbf{k+\kappa}}^{r}\ $with coordinates of $\mathbf{\kappa}$ are $0$ or
$1$ form an $r$-dimensional polygon with parallel faces, trapezoid if $r=2,$
denote it by $\Delta_{\mathbf{k}}^{r}.$

Now supplement the constructed splitting $\{\Delta_{\mathbf{k}}^{r}%
,\mathbf{k}=(k_{1},...k_{r},0,...,0)\in\mathbb{Z}^{d}\}$ of $M$ with splitting
of $M_{\varepsilon}.$ For any polygon $\Delta_{\mathbf{k}}^{r},$ consider the
hyperspace $\mathbf{t}_{\mathbf{k}}^{r}+\mathbb{R}^{d-r}$ and take in
$\mathbb{R}^{d-r}$ $d-r$-dimensional uniform rectangular net with sides all
equal to $\mathbf{q}_{\mathbf{t}_{\mathbf{k}}^{r}}(u)\mathbf{T}_{d-r}$, where
$\mathbf{T}_{d-r}=(0,...0,T,...,T)\in$ $\mathbb{R}^{d},$ with first $r$ zeros.
That is, denote $\mathbb{T}^{d}=[0,T]^{d},$ $\mathbb{T}^{d-r}=\mathbb{T}%
^{d}\cap\mathbb{R}^{d-r},$ and take in $\mathbf{t}_{\mathbf{k}}^{r}%
+\mathbb{R}^{d-r},$
\begin{equation}
\Delta_{\mathbf{kl}}^{d-r}=\mathbf{t}_{\mathbf{k}}^{r}+\mathbf{q}%
_{\mathbf{t}_{\mathbf{k}}^{r}}(u)\mathbf{lT}_{d-r}+\mathbf{q}_{\mathbf{t}%
_{\mathbf{k}}^{r}}(u)\mathbb{T}^{d-r},\ \mathbf{l}=(0,...,0,l_{r+1}%
,...,l_{d})\in\mathbb{Z}^{d}. \label{spitting}%
\end{equation}
Finaly denote
\begin{equation}
\Delta_{\mathbf{kl}}^{d}:=\Delta_{\mathbf{k}}^{r}\times\Delta_{\mathbf{kl}%
}^{d-r},\ \ \Delta_{\mathbf{kl}}^{d}\cap M_{\varepsilon}\neq\varnothing.
\label{splittingF}%
\end{equation}
Denoting $\mathbf{t}_{\mathbf{k}}^{d-r}:=\mathbf{q}_{\mathbf{t}_{\mathbf{k}%
}^{r}}(u)\mathbf{lT}_{d-r},$ since $\mathbb{T}^{d}=\mathbf{T}_{d}%
\otimes\mathbf{T}_{d-r},$ we can write this as
\[
\Delta_{\mathbf{kl}}^{d}=\mathbf{t}_{\mathbf{k}}^{r}+\mathbf{q}_{\mathbf{t}%
_{\mathbf{k}}}(u)((\mathbf{k+l)}\mathbb{T}+\mathbb{T}),\
\]
so that Lemma \ref{lemma_pickands} is applicable with $\mathbf{t=t}%
_{\mathbf{k}}^{r}$ and $T=$ $(\mathbf{k+l)}\mathbb{T}+\mathbb{T}.$ In case
$\mathbf{t}_{\mathbf{k}}^{r}\notin M$ (some end point of the net
$\mathcal{N}_{r}$), take $\mathbf{t}\in M\mathbf{,}$ another vertex of the
same polygon and let $\mathbf{q}_{\mathbf{t}_{\mathbf{k}}^{r}}(u)=\mathbf{q}%
_{\mathbf{t}}(u).$

Introduce the events%
\[
A_{\mathbf{k,l}}:=\{\sup_{\mathbf{t\in}\Delta_{\mathbf{k,l}}^{d}}%
X(\mathbf{t)}>u\}.
\]
We have for any $A\subset M_{\varepsilon},$%
\begin{equation}%
%TCIMACRO{\dbigcup \limits_{\Delta_{\mathbf{k,l}}\subset A}}%
%BeginExpansion
{\displaystyle\bigcup\limits_{\Delta_{\mathbf{k,l}}\subset A}}
%EndExpansion
A_{\mathbf{k,l}}\subseteq\left\{  \max_{\mathbf{s\in}A}X(\mathbf{s)>}%
u\right\}  \subseteq%
%TCIMACRO{\dbigcup \limits_{\Delta_{\mathbf{k,l}}\cap A\neq\varnothing}}%
%BeginExpansion
{\displaystyle\bigcup\limits_{\Delta_{\mathbf{k,l}}\cap A\neq\varnothing}}
%EndExpansion
A_{\mathbf{k,l}}. \label{bonferr}%
\end{equation}

\subsection{Locally homogeneous Gaussian fields.}

Denote by $X_{0}(\mathbf{t})$, a Gaussian a.s. continuous zeromean field with
variance $1$ and covariation function $r(\mathbf{s}_{1},\mathbf{s}_{2}).$ That
is,
\[
X(\mathbf{t})\overset{d}{=}\sigma(\mathbf{t})X_{0}(\mathbf{t}),
\]
in distributions. Assuming that Condition \ref{CondA4} is fulfilled for all
$\mathbf{t\in}\mathcal{S}_{1}\subset\mathcal{S},$ we can say that the field
$\xi(\mathbf{t}),$ $\mathbf{t\in}\mathcal{S}_{1},$ is \emph{locally
homogeneous. }In particular, in Condition \ref{CondA4}, \emph{ } the Gaussian
field $X(\mathbf{t}),$ $\mathbf{t}\in\mathcal{M\cap U},$ is locally
homogeneous. The following theorem is formulated for $X_{0}(\mathbf{t}),$ but
in fact it is valid for any locally homogeneous Gaussian fields satisfying the
above conditions.

\begin{theorem}
\label{Pickands Theorem} Let Condition \ref{CondA4} be fulfilled for all
$\mathbf{t\in}\mathcal{U}.$ Then for any $\mathcal{U}_{1}\mathcal{\subset U}$
which is the closure of an open set,
\begin{equation}
P(\max_{\mathbf{t\in}\mathcal{U}_{1}}X_{0}(\mathbf{t)>}u)=H_{\mathbf{q}}%
\int_{\mathcal{U}_{1}}\prod_{i=1}^{r}C_{\mathbf{t}i}^{-1}(u)d\mathbf{t}%
\prod_{i=1}^{r}q_{i}(u)^{-1}\Psi(u)(1+o(1)) \label{Picands0}%
\end{equation}
as $u\rightarrow\infty,$ with
\[
H_{\mathbf{q}}=\lim_{T\rightarrow\infty}T^{-d}H_{\mathbf{q}}\left(
[0,T]^{d}\right)  \in(0,\infty)\mathbf{.}%
\]
This assertion holds even if corresponding $U_{1}$ depends of $u$,
$U_{1}=U_{1}(u),$ provided there exist boxes $U_{1}^{\pm}(u)=\otimes_{i=1}%
^{d}[-U_{1i}^{\pm}(u),U_{1i}^{\pm}(u)]$ such that $U_{1}^{-}(u)\subset
U_{1}(u)\subset U_{1}^{+}(u)$ with $U_{1i}^{-}(u)q_{i}(u)\rightarrow\infty,$
$u\rightarrow\infty,$ $i=1,...,d.$
\end{theorem}

\begin{remark}
Remark that the assertion is valid also in case the parametric set
$\mathcal{U}_{1}$ expands unboudedly, with some restrictions on the expanding,
it follows from the proof. See corresponding conditions in \cite{KP}.
\end{remark}

\textbf{Outline of proof. } Introduce a positive increasing function
$\kappa(u),$ $u>0,$ with
\[
\lim_{u\rightarrow\infty}\kappa(u)=\infty,\ \text{but }\lim_{u\rightarrow
\infty}q_{i}(u)\kappa(u)=0,\ i=1,...,d.
\]
Cover the set $\mathcal{U}_{1}$ with boxes
\[
B_{\mathbf{k}}(u)=\kappa(u)\mathbf{kq}(u)+B_{\mathbf{0}}(u)\mathbf{,}%
\quad\mathbf{k}\in\mathbb{Z}^{d},u>0,
\]
where
\[
B_{\mathbf{0}}:=\kappa(u)\bigotimes_{i=1}^{d}[0,q_{i}(u)].
\]
Denote%
\[
\mathbf{C}_{\mathbf{k}}^{\pm}:=(C_{i\mathbf{k}}^{\pm}\mathbf{,}%
i=1,...,d),\mathbf{\ }%
\]
with
\[
C_{i\mathbf{k}}^{+}:=\max_{\mathbf{t\in}B_{\mathbf{k}}(u)}C_{i\mathbf{t}%
}\mathbf{,\ }C_{i}^{-}(\mathbf{k)}:=\min_{\mathbf{t\in}B_{\mathbf{k}}%
(u)}C_{i\mathbf{t}},\ i=1,...,d.\mathbf{\ }%
\]
Fix $\mathbf{k,}$ and consider Gaussian zeromean homogeneous fields
$\xi_{\mathbf{k}}^{\pm}(\mathbf{t}),$ $\mathbf{t\in}B_{\mathbf{k}}(u),$ with
covariance functions $r_{\mathbf{k}}^{\pm}(\mathbf{t})$ satisfying Condition
\ref{CondA4} with $\mathbf{C}_{\mathbf{k}}^{\pm}$ instead of $\mathbf{C}%
_{\mathbf{t}},$ correspondingly. Applying the constructed above splitting of
$\ B_{\mathbf{k}}(u)$ with polygons $\Delta_{\mathbf{i}}^{d}$ and using Lemma
\ref{lemma_pickands}, we get similarly to the Theorem 7.1, \cite{book}, proof,
the asymptotic behavior of
\[
P(\max_{\mathbf{t\in}B_{\mathbf{k}}(u)}\xi_{\mathbf{k}}^{\pm}(\mathbf{t}%
)\mathbf{>}u).
\]
Then, using Slepian inequality to bound the probability for $X_{0}%
(\mathbf{t),}$\ $\mathbf{t\in}B_{\mathbf{k}}(u),$ and applying again Double
Sum Method for the probability in question, we get Theorem. See quite similar
proofs in \cite{KP}, \cite{KHP}, as well as in \cite{book}, \cite{lectures}.

Remark that since (\ref{PUu}), it is easier to perform the proof for the
corresponding $U_{1}=\varphi^{-1}(\mathcal{U}_{1}),$ which leads to the
integral over $U_{1}$ with the functions $C_{\mathbf{t}i}$ and $q_{i}(u)$
constructed for $U=\varphi^{-1}(\mathcal{U}_{1})$ and the fixed at the
agreement coordinate axes. Finally we simply pass to the integral over
$\mathcal{M}_{1},$ with another, of course, volume element $d\mathbf{t}$
(changing variables $\mathbf{t}\rightarrow\varphi(\mathbf{t})$), \emph{but
leaving the same} $\mathbf{q}(u)$ \emph{constructed in Condition} \ref{CondA4}.

We need also the following corollary of Theorem \ref{Pickands Theorem}.

\begin{proposition}
\label{Pickands} In the assumptions of this section, for any $\mathcal{M}%
_{1}\mathcal{\subset U\cap M}$ which is the closure of an open set,
\begin{equation}
P(\max_{\mathbf{t\in}\mathcal{M}_{1}}X(\mathbf{t)>}u)=H_{\mathbf{q,}r}%
\int_{\mathcal{M}_{1}}\prod_{i=1}^{r}C_{\mathbf{t}i}^{-1}(u)d\mathbf{t}%
\prod_{i=1}^{r}q_{i}(u)^{-1}\Psi(u)(1+o(1)) \label{Picands1}%
\end{equation}
as $u\rightarrow\infty,$ with
\[
H_{\mathbf{q,}r}=\lim_{T\rightarrow\infty}T^{-r}H_{\mathbf{q}}\left(
[0,T]^{r}\times\lbrack0,0]^{d-r}\right)  \in(0,\infty)\mathbf{.}%
\]
This assertion holds even if corresponding $M_{1}$ depends of $u$,
$M_{1}=M_{1}(u),$ like in the Theorem \ref{Pickands Theorem}.
\end{proposition}

\subsection{Homogeneous like Gaussian fields.}

Consider the homogeneous like case, that is $h_{1\mathbf{t}}(\mathbf{s)=}0$
for all $\mathbf{s}\in M_{\varepsilon}$ and $\mathbf{t}\in M\mathbf{.}$ Take
$\varepsilon>0$ sufficiently small, and denote
\begin{equation}
f(\mathbf{s})=\frac{1}{2}(1-\sigma^{2}(\mathbf{s})),\ \mathbf{s}\in
M_{\varepsilon}. \label{f(t)}%
\end{equation}
Introduce the Laplace type integral,%
\begin{equation}
L_{f}(\lambda):=\int_{U}e^{-\lambda f(\mathbf{s)}}d\mathbf{s},\quad\lambda>0.
\label{laplace}%
\end{equation}
Remark that by standard asymptotic Laplace method, its asymptotic behavior as
$\lambda\rightarrow\infty$ depends only on behavior $f(\mathbf{t})$ at zero,
that is, as $\rho(\mathbf{s},M)\rightarrow0.$

\begin{proposition}
\label{stationary like} Let Conditions \ref{CondA1}--\ref{CondA5} be
fulfilled. If further $h_{1}(\mathbf{t})=0$ for all $\mathbf{t}$, that is,
$\cup_{\mathbf{t}\in M}K_{0\mathbf{t}}=U,$ we have for any map $\mathcal{U},$
from the fixed atlas on $\mathcal{M}_{\varepsilon},$
\begin{equation}
P(\mathcal{U};u)=(1+o(1))H_{\mathbf{q}}\int_{\mathcal{M\cap U}}\prod_{i=1}%
^{d}C_{\mathbf{t}i}^{-1}(u)d\mathbf{t}L_{f}(u^{2})\ \prod_{i=1}^{d}q_{i}%
^{-1}(u)\Psi(u)\ \ \label{case S}%
\end{equation}
as $u\rightarrow\infty$.
\end{proposition}

Proof of this proposition is based on the same constructions as in Theorem
\ref{Pickands Theorem} proof. The only for the approximation of
$P(B_{\mathbf{k}}(u);u)$ the probabilities $P(\max_{\mathbf{t\in}%
B_{\mathbf{k}}(u)}\sigma_{\mathbf{k}}^{\pm}X_{0}(\mathbf{t})\mathbf{>}u)$ are
used directly, with $\sigma_{\mathbf{k}}^{\pm},$ upper and bound bounds of
$\sigma(\mathbf{t}),$ $\mathbf{t\in}B_{\mathbf{k}}(u)$, with followed Theorem
\ref{Pickands Theorem} application and standard limit passage to integral as
$u\rightarrow\infty.$ Very similar detailed proof in case $\dim\mathcal{M}=0$
is given in \cite{KP}, Proposition 3, which combined both above proofs for
locally homogeneous and homogeneous like cases. Notice that the proof is
performed for the manifold $\mathcal{M}_{u},$ $u\rightarrow\infty$ with
following application Lemma \ref{extracting}. Most technical detailes here are
related to double sums estimations in corresponding Bonferroni inequalities.

\subsection{Talagrand case.}

\begin{proposition}
\label{Talagrand_case}In Conditions \ref{CondA1}--\ref{CondA5}, if
$\cup_{\mathbf{t}\in M}K_{\infty\mathbf{t}}=M_{\varepsilon},$ then for any map
$\mathcal{U}$ from the fixed atlas on $\mathcal{M}_{\varepsilon},$
\begin{equation}
P(\mathcal{U};u)=H_{\mathbf{q,}r}\int_{\mathcal{M\cap U}}\prod_{i=1}%
^{r}C_{\mathbf{t}i}^{-1}(u)d\mathbf{t}\prod_{i=1}^{r}q_{i}(u)^{-1}%
\Psi(u)(1+o(1)),\ \ \ \label{case T}%
\end{equation}
as $u\rightarrow\infty.$ The constant $H_{\mathbf{q,}r}$ is defined in
(\ref{Picands1}).
\end{proposition}

In other words, the asymptotical behavior of the probability in question is
coincides with the asymptotic behavior of $P(M;u),$ compare with
(\ref{Picands1}). Notice that it is proved in \cite{KP} that in case $M$
consists of one point the right part is equal simply to $\Psi(u)(1+o(1)).$
Therefore it follows from Lemma \ref{extracting} and trivial application of
Bonferroni inequalities that in case $\dim M=0$ since $M$ is finitely
connected, the probability is equivalent to $\operatorname*{card}(M)\Psi(u),$
$u\rightarrow\infty.$ See Remark \ref{r=0}.

\textbf{Outline of proof. }In this case, the informative set $M_{u}$ is
contained in the two closest to $M$ layers of the constructed above splitting
and moreover of any small thickness. Namely, take the same splitting
$\{\Delta_{\mathbf{k}}^{r},\mathbf{k}=(k_{1},...k_{r},0,...,0)\in
\mathbb{Z}^{d}\}$ of $M$, but now denote $\mathbf{\epsilon}_{d-r}%
=(0,...0,\varepsilon,...,\varepsilon)\in$ $\mathbb{R}^{d},$ $\varepsilon>0,$
with first $r$ zeros. Similarly the abobe, denote $\mathbb{E}^{d}%
=[-\varepsilon,\varepsilon]^{d},$ $\mathbb{E}^{d-r}=\mathbb{E}^{d}%
\cap\mathbb{R}^{d-r},$ and take in $\mathbf{t}_{\mathbf{k}}^{r}+\mathbb{R}%
^{d-r},$
\[
\Delta_{\mathbf{k0}}^{d-r}=\mathbf{t}_{\mathbf{k}}^{r}+\mathbf{q}%
_{\mathbf{t}_{\mathbf{k}}^{r}}(u)\mathbb{E}^{d-r},
\]
compare with (\ref{spitting}). Finally denote again
\[
\Delta_{\mathbf{k0}}^{d}:=\Delta_{\mathbf{k}}^{r}\times\Delta_{\mathbf{k0}%
}^{d-r},\ \ \Delta_{\mathbf{k0}}^{d}\cap M_{u}\neq\varnothing,
\]
and%
\[
A_{\mathbf{k,0}}:=\{\sup_{\mathbf{t\in}\Delta_{\mathbf{k,0}}^{d}}%
X(\mathbf{t)}>u\},A_{\mathbf{k}}^{r}:=\{\sup_{\mathbf{t\in}\Delta_{\mathbf{k}%
}^{r}}X(\mathbf{t)}>u\}.\
\]
We have in Talagrand case, that for any $A\subset U\cap M_{u},$ any small
$\varepsilon>0,$ and all sufficiently large $u,$%
\[%
%TCIMACRO{\dbigcup \limits_{\Delta_{\mathbf{k}}^{r}\subset M\cap A}}%
%BeginExpansion
{\displaystyle\bigcup\limits_{\Delta_{\mathbf{k}}^{r}\subset M\cap A}}
%EndExpansion
A_{\mathbf{k}}^{r}\subseteq\left\{  \max_{\mathbf{s\in}A}X(\mathbf{s)>}%
u\right\}  \subseteq%
%TCIMACRO{\dbigcup \limits_{\Delta_{\mathbf{k,0}}\cap A\neq\varnothing}}%
%BeginExpansion
{\displaystyle\bigcup\limits_{\Delta_{\mathbf{k,0}}\cap A\neq\varnothing}}
%EndExpansion
A_{\mathbf{k,0}}.
\]
Asymptotic behavior of the left side event probability have been already
computed, Proposition \ref{Pickands}. The right hand probability is bounded by
the sum of probabilities of $A_{\mathbf{k,0}}$ with following Lemma
\ref{lemma_pickands} application and then taking $\varepsilon$ arbitrarily
small. This part is absolutely similar to the corresponding part of Theorem
8.2, \cite{book} (Lemma 8.4), see also \cite{lectures}, Lecture 10, and
\cite{KP}.

\subsection{Transition case.}

In this case we have to consider separately two parts of $M,$ the first one
belongs to $\partial U$, and the second one located inside of $U.$ Denote
$\Pi_{\mathbf{t}}:=[0,T]^{r}\times\lbrack-S,S]^{d-r}\mathbf{I}%
_{\mathbf{t\notin}\partial U}+[0,T]^{r}\times\lbrack0,S]^{d-r}\mathbf{I}%
_{\mathbf{t\in}\partial U},$ $\mathbf{t}\in M,$
\[
P_{\mathbf{q,t}}(\Pi\mathbf{)=}E\exp(\max_{\mathbf{s\in}\Pi_{\mathbf{t}}}%
\chi(\mathbf{s})-h_{1\mathbf{t}}(\mathbf{s)}).\mathbf{\ \ }%
\]

\begin{proposition}
\label{PPcase}In the above conditions, if $\cup_{\mathbf{t}\in M}%
K_{c\mathbf{t}}=M_{\varepsilon},$ then%
\begin{equation}
P(U;u)=\int_{\mathcal{M}}P_{\mathbf{q,t}}\prod_{i=1}^{r}C_{\mathbf{t}i}%
^{-1}(u)d\mathbf{t}\prod_{i=1}^{r}q_{i}(u)^{-1}\Psi(u)(1+o(1)) \label{case P}%
\end{equation}
as $u\rightarrow\infty,$ with
\[
P_{\mathbf{q,t}}^{\nu}=\lim_{T,S\rightarrow\infty}T^{-r}P_{\mathbf{q,t}}^{\nu
}(\Pi_{\mathbf{t}}\mathbf{)}\in(0,\infty),
\]
see (\ref{pit}).
\end{proposition}

\textbf{Outline of proof. }We again concentrate ourselves on the construction
of corresponding splitting. In this case the informative set $M_{u}$ is much
wider than in the previous case. Take again the same splitting $\{\Delta
_{\mathbf{k}}^{r},\mathbf{k}=(k_{1},...k_{r},0,...,0)\in\mathbb{Z}^{d}\}$ of
$M$, but now denote $\mathbf{S}_{d-r}=(0,...0,S,...,S)\in$ $\mathbb{R}^{d},$
$S>0,$ with first $r$ zeros, $S$ will be taken large. Denote $\mathbb{S}%
^{d}=[-S,S]^{d},$ $\mathbb{S}^{d-r}=\mathbb{S}^{d}\cap\mathbb{R}^{d-r},$ and
repeat (\ref{spitting}, \ref{splittingF}), taking
\[
\Delta_{\mathbf{kl}}^{d-r}=\mathbf{t}_{\mathbf{k}}^{r}+\mathbf{q}%
_{\mathbf{t}_{\mathbf{k}}^{r}}(u)\mathbf{lS}_{d-r}+\mathbf{q}_{\mathbf{t}%
_{\mathbf{k}}^{r}}(u)\mathbb{S}^{d-r},\mathbf{l}=(0,...,0,l_{d-r+1}%
,...,l_{d})\in\mathbb{Z}^{d}%
\]
and
\[
\Delta_{\mathbf{kl}}^{d}:=\Delta_{\mathbf{k}}^{r}\times\Delta_{\mathbf{kl}%
}^{d-r},\ \ \Delta_{\mathbf{kl}}^{d}\cap M_{u}\neq\varnothing.
\]
Denote again%
\[
A_{\mathbf{k,l}}:=\{\sup_{\mathbf{t\in}\Delta_{\mathbf{k,l}}^{d}}%
X(\mathbf{t)}>u\},\
\]
We have in this case, that for any $A\subset M_{\varepsilon},$ and all
sufficiently large $u,$%
\[
\bigcup\limits_{\Delta_{\mathbf{k0}}^{d}\subset A}A_{\mathbf{k0}}^{r}%
\subseteq\left\{  \max_{\mathbf{s\in}A}X(\mathbf{s)>}u\right\}  \subseteq
\bigcup\limits_{\Delta_{\mathbf{k0}}^{d}\cap A\neq\varnothing}A_{\mathbf{k,0}%
}\ \cup\bigcup\limits_{\Delta_{\mathbf{kl}}^{d}\cap A\neq\varnothing
,\mathbf{l\neq0}}A_{\mathbf{k,l}}.
\]
Asymptotic behavior of the probability of left side event can be evaluated by
standard Double Sum Method, using Lemma \ref{lemma_pickands} and standard
estimation of the corresponding double sum, which is infinitely smaller that
the single one as $u\rightarrow\infty\ $for non-closed $\Delta_{\mathbf{k0}%
}^{d}$ and the rest one, with closed $\Delta_{\mathbf{k0}}^{d},$ becames small
as $T\rightarrow\infty,$ even with fixed $S,$ the same for the first summand
in the right hand part. See \cite{book}, \cite{lectures}. The probability of
the second summand in the right hand part can be estimated from above by sum
of probabilities of $A_{\mathbf{k,l}}.$Then one should estimate theirs
asymptotic behavior as $\allowbreak u\rightarrow\infty,$ and then it turns out
that this estimation is negligibly smaller as $S\rightarrow\infty.$This
estimation is practically identical to the correponding place in \cite{KP}.
See again the corresponding parts of \cite{book}, \cite{lectures}.

\subsection{General case.}

Denote
\begin{align*}
K_{0}(M)  &  =\{\mathbf{t}\in M:K_{0\mathbf{t}}\neq\varnothing\},K_{c}%
(M)=\{\mathbf{t}\in M:K_{0\mathbf{t}}=\varnothing,K_{c\mathbf{t}}%
\neq\varnothing\},\\
K_{\infty}(M)  &  =\{\mathbf{t}\in M:K_{0\mathbf{t}}=\varnothing
,K_{c\mathbf{t}}=\varnothing,K_{\infty\mathbf{t}}\neq\varnothing\},
\end{align*}
see (\ref{cases}). Furthermore, denote
\[
K_{0}:=%
%TCIMACRO{\dbigcup \limits_{\mathbf{t}\in K_{0}(M)}}%
%BeginExpansion
{\displaystyle\bigcup\limits_{\mathbf{t}\in K_{0}(M)}}
%EndExpansion
K_{0\mathbf{t}},K_{c}:=%
%TCIMACRO{\dbigcup \limits_{\mathbf{t}\in K_{c}(M)}}%
%BeginExpansion
{\displaystyle\bigcup\limits_{\mathbf{t}\in K_{c}(M)}}
%EndExpansion
K_{c\mathbf{t}},K_{\infty}:=%
%TCIMACRO{\dbigcup \limits_{\mathbf{t}\in K_{\infty}(M)}}%
%BeginExpansion
{\displaystyle\bigcup\limits_{\mathbf{t}\in K_{\infty}(M)}}
%EndExpansion
K_{\infty\mathbf{t}}.
\]

Assume that all these sets are simply designed in the following sense.

\begin{condition}
\label{CondA6} \textbf{(Structure Condition)} Assume that all these sets are
finitely connected smooth (two times differentiable) manifolds. Furthermore,
assume that $K_{0}$ is either empty or consists of finite number of smooth
disjoint manifolds with positive dimensions, namely,
\begin{equation}
K_{0}=\bigcup_{i=1}^{n}\mathcal{K}_{0i},\ \dim\mathcal{K}_{0i}=k_{i}%
,\ 0<k_{1}\leq k_{2}\leq...\leq k_{n}\leq d. \label{part}%
\end{equation}
Assume also that for any $i,$ $k_{i}$-dimensional volume of $\mathcal{K}_{0i}$
is finite, $|\mathcal{K}_{0i}|<\infty.$
\end{condition}

We use here caligrafic $\mathcal{K}_{\cdot}$ because of these manifolds in $U$
can be indeed of general type, not necessary linear subspaces. Fix $i$ with
$k_{i}<d,$ and consider in $\mathcal{K}_{0i}$ curvilinear coordinates. For
$\mathbf{t\in}\mathcal{K}_{0i},$ using Proposition \ref{e_reg_var}, choose
coordinate vectors $\mathbf{e}_{j}(\mathbf{t),}$ $j=1,...,k_{i}\mathbf{\ }$of
this curvilinear coordinates and complete them to a basis $\{\mathbf{e}%
_{j}(\mathbf{t),}$ $j=1,...,k_{i},$ $\mathbf{\tilde{e}}_{j}(\mathbf{t),}$
$j=k_{i}+1,...,d\}$ in $\mathbb{R}^{d}.$ Write vector $\mathbf{q}_{\mathbf{t}%
}(u)$ in these coordinates, remind that $\mathbf{q}_{\mathbf{t}}%
(u)=\mathbf{C}_{\mathbf{t}}\mathbf{q}_{\mathbf{0}}(u),$ see \ref{CondA4},
\begin{equation}
\mathbf{q}_{\mathbf{t}}^{i}(u)=(C_{\mathbf{t}j}^{i}q_{\mathbf{0}j}%
^{i}(u),\ j=1,...,k_{i},C_{\mathbf{t}j}^{i}q_{\mathbf{0}j}^{i}(u),j=k_{i}%
+1,...,d), \label{K01}%
\end{equation}
where, in accordance with Proposition \ref{e_reg_var}, $(q_{\mathbf{0}j}%
^{i},j=1,...,d)$ is just an permutation of $(q_{\mathbf{0}j},j=1,...,d),$
whereas $C_{\mathbf{t}j}^{i}$ are other positive continuous bounded functions.
Denote corresponding positive limits by
\begin{equation}
h_{\mathbf{t}}^{i}(\mathbf{s}):=\lim_{u\rightarrow\infty}u^{2}(1-r(\mathbf{q}%
_{\mathbf{t}}^{i}(u)\mathbf{s}), \label{K02}%
\end{equation}
where $\mathbf{s}=(s_{1},...,s_{d})$ is written in these coordinates. In fact
$h_{\mathbf{t}}^{i}(\mathbf{s})=h(O_{i,\mathbf{t}}\mathbf{s),}$ where
$O_{i,\mathbf{t}}$ is an orthogonal transition matrix to the curvilinear
coordinates with the orthogonal complement to a basis in $\mathbb{R}^{d},$ as above.

We have,%
\begin{equation}
\lim_{u\rightarrow\infty}\frac{1-\sigma(\mathbf{q}_{\mathbf{t}}^{i}%
(u)\mathbf{s})}{1-r(\mathbf{q}_{\mathbf{t}}^{i}(u)\mathbf{s})}=0. \label{K03}%
\end{equation}
By analogy with proofs of Theorem \ref{Pickands Theorem} and Proposition
\ref{stationary like} proof, we build a partition of $\mathcal{K}_{0i}$ with
$k_{i}$-dimensional blocks, similar to $B_{0}$ and $B_{k},$ with the same
$\kappa(u)$ and $\mathbf{q}_{\mathbf{t}}^{i}(u)$ from (\ref{K01}). Denote
these blocks $B(\mathbf{t}_{\nu}),$ where $\{\mathbf{t}_{\nu},\nu=1,...,N\}$
is the corresponding to this partition grid, that is,
\begin{equation}
\mathcal{K}_{0i}=\bigcup_{\nu=1}^{N}B(\mathbf{t}_{\nu}). \label{grid}%
\end{equation}
Similarly to the same proofs, but for $\mathcal{K}_{0i}\cap U\cap M_{u}$
instead of $U\cap M_{u},$ using Theorem \ref{Pickands Theorem}, the second
part we get asymptotic behavior of the probability $P(B(\mathbf{t}_{\nu});u)$
for all $\nu.$ Then, thinning the grid unboundedly, we get, using Condition
\ref{CondA6}, the following Lemma.

\begin{lemma}
\label{K0i} For any $\mathcal{K}_{0i}$ from the partition (\ref{part}) with
$k_{i}<d,$
\[
P(\mathcal{K}_{0i}\cap U;u)=\int_{\mathcal{K}_{0i}\cap U}H_{\mathbf{q}%
_{\mathbf{t}}^{i}}\prod_{j=1}^{k_{i}}(q_{\mathbf{t}j}^{i}(u)^{-1}%
e^{-u^{2}f(\mathbf{t)}}d\mathbf{t}\Psi(u)(1+o(1))\ \
\]
as $u\rightarrow\infty$. Remind that, by our agreement, here $d\mathbf{t}$ is
an elementary $k_{i}$-dimensional volume of $\mathcal{K}_{0i},$ the
corresponding integration domain, and $f(\mathbf{t})$ is given by (\ref{f(t)}).
\end{lemma}

Turning to (\ref{part}), we have, using again Bonferroni inequalities, with
some additional standard techniques for probabilities of large excursions, see
\cite{KP}, \cite{KPR}, we get the following.

\begin{proposition}
\label{Th_dimK_0<d} If $\dim\mathcal{K}_{0i}<d,$ $i=1,...,n,$
\begin{equation}
P(K_{0}\cap U;u)=(1+o(1))\sum_{i=1}^{n}P(\mathcal{K}_{0i}\cap
U;u),\ \ \label{K0_fin2}%
\end{equation}
as $u\rightarrow\infty$.
\end{proposition}

\begin{remark}
\label{dimK d} Remark that if for some $i$ in (\ref{part}), $\dim
\mathcal{K}_{0i}=d,$ Proposition \ref{stationary like} should be used with
another Laplace integral (\ref{laplace}), namely,
\begin{equation}
L_{f}(\lambda):=\int_{\mathcal{K}_{0}}e^{-\lambda f(\mathbf{t)}}%
d\mathbf{t},\quad\lambda>0. \label{laplace1}%
\end{equation}

\end{remark}

\begin{remark}
\label{dim K d}Remark that the summands in (\ref{K0_fin2}) can have different
orders in $u$ depending on the dimension of the corresponding component
$\mathcal{K}_{0i},$ on behavior of $q_{j}^{i,\mathbf{t}}(u)$s and on behavior
of $\sigma^{2}(\mathbf{t})$. Hence only summands with slowest order play a
role. Remark also that by Remark \ref{dimK d}, if for some $i,$ $\dim
\mathcal{K}_{0i}=d,$ no summands with $\dim\mathcal{K}_{0i}<d$ in
(\ref{K0_fin2}) contribute to the asymptotic behavior of $P(U;u).$
\end{remark}

\begin{remark}
\label{other cases} Evaluations of asymptotic behaviors of $P(\mathcal{K}%
_{\infty},u)$ and $P(\mathcal{K}_{c},u)$ are obvious, with using Condition
\ref{CondA6} and Propositions \ref{Talagrand_case} and \ref{PPcase}, respectively.
\end{remark}

\section{Conclusion. Summary of main results.}

Thus, using again Bonferroni, we get for any map $\mathcal{U}$ from the atlas
$\mathcal{A}$ on $\mathcal{M}_{\varepsilon}$ that
\[
P(\mathcal{U},u)=P(U,u)=
\]%
\[
=(P(K_{0}\cap U;u)+P(K_{c}\cap U;u)+P(K_{\infty}\cap U;u))(1+o(1))
\]
as $u\rightarrow\infty$, and the asymptotic behavior of the first summand is
evaluated in Proposition \ref{Th_dimK_0<d}, for the asymptotic behaviors of
the rest two summands see Remark \ref{other cases}.

The final steps are

1) to scale the atlas making maps pairwise nonintersected, namely number the
maps, $\mathcal{U}_{i},$ $i=1,...,M,$ $M<\infty,$ and take $\mathcal{U}_{1},$
$\mathcal{U}_{2}^{\prime}=\mathcal{U}_{2}\setminus\mathcal{U}_{1},$ $...,$
$\mathcal{U}_{k}^{\prime}=\mathcal{U}_{k}\setminus(\cup_{i=1}^{k-1}%
\mathcal{U}_{i}^{\prime}),...$ as the new atlas;

2) to wtite
\[
\{\max_{\mathbf{t\in}\mathcal{M}_{\varepsilon}}X(\mathbf{t})>u\}=\bigcup
\limits_{\mathcal{U}^{\prime}\in\mathcal{A}}\{\max_{\mathbf{t\in}%
\mathcal{U}^{\prime}}X(\mathbf{t})>u\};
\]

3) to use Bonferroni inequalities with estimating of the double sum of the
joint probabilities;

and

4) to use back passage from the informative manifold $\mathcal{M}%
_{\varepsilon}$ to whole $\mathcal{S}.$

Remark that our results here agree with assertions of Theorem 3, \cite{KHP},
where $d=1,$ and Theorem 1 in \cite{KP}, \cite{KPR} for $\dim\mathcal{M}=0.$

\bigskip

\end{document}